\documentclass{article}

\usepackage{amsmath, amssymb}

\title{A note on the Grover walk and the generalized Ihara zeta function of the one-dimensional integer lattice} 

\author{Takashi KOMATSU \\
Department of Bioengineering School of Engineering,\\ 
The University of Tokyo \\
Bunkyo, Tokyo, 113-8656, JAPAN \\ 
e-mail: komatsu@coi.t.u-tokyo.ac.jp \\ 
Norio KONNO \\
Department of Applied Mathematics, Faculty of Engineering, \\ 
Yokohama National University \\
Hodogaya, Yokohama 240-8501, JAPAN \\
e-mail: konno-norio-bt@ynu.ac.jp, Tel.: +81-45-339-4205, Fax: +81-45-339-4205 \\ 
Iwao SATO \\ 
Oyama National College of Technology \\
Oyama, Tochigi 323-0806, JAPAN \\ 
e-mail: isato@oyama-ct.ac.jp }
 
\begin{document}
 \maketitle

\clearpage

\vspace{5mm}



\vspace{5mm}

The contact author for correspondence: 

Iwao Sato 

Oyama National College of Technology, 
Oyama, Tochigi 323-0806, JAPAN

Tel: +81-285-20-2176

Fax: +81-285-20-2880

E-mail: isato@oyama-ct.ac.jp

\clearpage

\begin{abstract} 
Chinta, Jorgenson and Karlsson introduced a generalized version of the determinant
formula for the Ihara zeta function associated to finite or infinite regular graphs. On
the other hand, Konno and Sato obtained a formula of the characteristic polynomial
of the Grover matrix by using the determinant expression for the second weighted zeta
function of a finite graph. In this paper, we focus on a relationship between the Grover
walk and the generalized Ihara zeta function. That is to say, we treat the generalized
Ihara zeta function of the one-dimensional integer lattice as a limit of the Ihara zeta
function of the cycle graph.
\end{abstract}

\section{Introduction}

Ihara \cite{Ihara} introduced the Ihara zeta functions of graphs, and showed that the reciprocal of 
the Ihara zeta function of a regular graph is an explicit polynomial. 
Afterwards, the Ihara zeta function of a finite graph was studied in \cite{Serre, Sunada, Sunada2, Hashimoto, Bass, EHSW}.  
Furthermore, the Ihara zeta function of a finite graph was extended to an infinite graph 
in \cite{Bass, CM, GZ, GIL, GIL2}, and its determinant expressions were presented. 
Chinta, Jorgenson and Karlsson \cite{CJK} gave a generalized version of 
the determinant formula for the Ihara zeta function associated to 
finite or infinite regular graphs.

A discrete-time quantum walk is a quantum analog of the classical random walk on a graph 
whose state vector is governed by a matrix called the time evolution matrix. 
The time evolution matrix of a discrete-time quantum walk in a graph 
is closely related to the Ihara zeta function of a graph. 
Ren et al. \cite{RAEWH} gave a relationship between the discrete-time quantum walk and the Ihara zeta function of a graph.  
Konno and Sato \cite{KS} obtained a formula of the characteristic polynomial of the Grover matrix 
by using the determinant expression for the second weighted zeta function of a graph. 

In this paper, we consider the relation between the Grover walk and the generalized
Ihara zeta function, and present the generalized Ihara zeta function of the one-dimensional
integer lattice as a limit of the Ihara zeta function of the cycle graph. 

In Section 2, we state a review for the Ihara zeta function of a finite graph and the generalized Ihara zeta function 
of a finite or infinite vertex transitive graph. 
In Section 3, we deal with the Grover walk on a graph as a discrete-time quantum walk on a graph. 
In Section 4, we treat the generalized Ihara zeta function of $ \mathbb{Z} $ as a limit of the Ihara zeta function 
of the cycle graph $C_n $ with $n$ vertices.

\section{The Ihara zeta function of a graph}

All graphs in this paper are assumed to be simple. 
Let $G=(V(G),E(G))$ be a connected graph (without multiple edges and loops) 
with the set $V(G)$ of vertices and the set $E(G)$ of unoriented edges $uv$ 
joining two vertices $u$ and $v$.
For $uv \in E(G)$, an arc $(u,v)$ is the oriented edge from $u$ to $v$. 
Let $D_G$ be the symmetric digraph corresponding to $G$. 
Set $D(G)= \{ (u,v),(v,u) \mid uv \in E(G) \} $. 
For $e=(u,v) \in D(G)$, set $u=o(e)$ and $v=t(e)$. 
Furthermore, let $e^{-1}=(v,u)$ be the {\em inverse} of $e=(u,v)$. 
For $v \in V(G)$, the {\em degree} $\deg {}_G v = \deg v = d_v $ of $v$ is the number of vertices 
adjacent to $v$ in $G$.  

A {\em path $P$ of length $n$} in $G$ is a sequence 
$P=(e_1, \cdots ,e_n )$ of $n$ arcs such that $e_i \in D(G)$,
$t( e_i )=o( e_{i+1} )(1 \leq i \leq n-1)$. 
If $e_i =( v_{i-1} , v_i )$ for $i=1, \cdots , n$, then we write 
$P=(v_0, v_1, \cdots ,v_{n-1}, v_n )$. 
Set $ \mid P \mid =n$, $o(P)=o( e_1 )$ and $t(P)=t( e_n )$. 
Also, $P$ is called an {\em $(o(P),t(P))$-path}. 
We say that a path $P=( e_1 , \cdots , e_n )$ has a {\em backtracking} 
if $ e^{-1}_{i+1} =e_i $ for some $i(1 \leq i \leq n-1)$. 
A $(v, w)$-path is called a {\em $v$-cycle} 
(or {\em $v$-closed path}) if $v=w$. 
A cycle $C=( e_1 , \ldots , e_r )$ has a {\em tail} if $e_r = e^{-1}_1 $. 
A cycle $C$ is {\em reduced} if $C$ has a neither a backtracking nor a tail. 
For a natural number $k \in \mathbb{N} $, let $N_k $ be the number of reduced cycles of length $k$ in $G$. 

The {\em Ihara zeta function} of a graph $G$ is 
a function of a complex variable $u$ with $|u|$ 
sufficiently small, defined by 
\[
{\bf Z} (G, u)= {\bf Z}_G (u)= \exp ( \sum^{\infty}_{k=1} \frac{N_k }{k} u^k ) . 
\]

Let $G$ be a connected graph with $n$ vertices $v_1, \cdots ,v_n $. 
The {\em adjacency matrix} ${\bf A}= {\bf A} (G)=(a_{ij} )$ is 
the square matrix such that $a_{ij} =1$ if $v_i$ and $v_j$ are adjacent, 
and $a_{ij} =0$ otherwise.
The {\em degree} of a vertex $v_i$ of $G$ is defined by 
$ \deg v_i = \deg {}_G v_i = \mid \{ v_j \mid v_i v_j \in E(G) \} \mid $. 
If $ \deg {}_G v=k$(constant) for each $v \in V(G)$, then $G$ is called 
{\em $k$-regular}.

\newtheorem{theorem}{Theorem}
\begin{theorem}[Ihara; Bass] 
Let $G$ be a connected graph. 
Then the reciprocal of the Ihara zeta function of $G$ is given by 
\[
{\bf Z} (G,u )^{-1} =(1- u^2 )^{r-1} 
\det ( {\bf I} -u {\bf A} (G)+ u^2 ( {\bf D} - {\bf I} )) , 
\]
where $r$ is the Betti number of $G$, 
and ${\bf D} =( d_{ij} )$ is the diagonal matrix 
with $d_{ii} = \deg v_i$ and $d_{ij} =0, i \neq j , 
(V(G)= \{ v_1 , \cdots , v_n \} )$. 
\end{theorem}

Let $G=(V(G),E(G))$ be a vertex transitive $(q+1)$-regular graph and $ x_0 \in V(G)$ 
a fixed vertex.  
Then the {\em generalized Ihara zeta function} $\zeta {}_G (u)$ of $G$ is defined by 
\[
\zeta {}_G (u)= \exp ( \sum^{\infty}_{m=1} \frac{N^0_m }{m} u^m ) , 
\]
where $N^0_m $ is the number of reduced $x_0$-cycles of length $m$ in $G$. 
Note that, for a finite graph, the classical Ihara zeta function is just the above
Ihara zeta function raised to the power equaling the number of vertices. 
Furthermore, the {\em Laplacian} of $G$ is given by 
\[
\Delta = \Delta (G) = {\bf D} - {\bf A} (G). 
\]

A formula for the generalized Ihara zeta function of a vertex transitive graph is given as follows:

\begin{theorem}[Chinta, Jorgenson and Karlsson] 
Let $G$ be a  vertex transitive $(q+1)$-regular graph with spectral measure $\mu $ for the Laplacian. 
Then 
\[
\zeta {}_G (u)^{-1} =(1-u^2 )^{(q-1)/2} \exp ( \int \log (1-(q+1- \lambda )u+q u^2 ) d \mu ( \lambda )) . 
\]
\end{theorem}

\section{The Grover walk on a graph}

Let $G$ be a connected graph with $n$ vertices and $m$ edges. 
Set $V(G)= \{ v_1 , \ldots , v_n \} $ and $d_j = d_{v_j} = \deg v_j , \ j=1, \ldots , n$. 
For $u \in V(G)$, let $D(u)= \{ e \in D(G) \mid t(e)=u \} $. 
Furthermore, let $\alpha {}_u , \ u \in V(G)$ be a unit vector with respect to $D(u)$, 
that is, 
\[
\alpha {}_u (e) =\left\{
\begin{array}{ll}
non \ zero \ complex \ number \ & \mbox{if $e \in D(u)$, } \\
0 & \mbox{otherwise, }
\end{array}
\right.
\]
where $ \alpha {}_u (e) $ is the entry of $\alpha {}_u $ corresponding to the arc $e \in D(G)$. 

Now, a $2m \times 2m$ matrix ${\bf C} $ is given as follows: 
\[
{\bf C} =2 \sum_{u \in V(G)} | \alpha {}_u \rangle \langle \alpha {}_u | - {\bf I}_{2m } . 
\]
The matrix {\bf C} is the {\em coin  operator} of the considered quantum walk. 
Note that ${\bf C} $ is unitary. 
Then the {\em time evolution matrix} ${\bf U} $ is defined by 
\[
{\bf U} = {\bf S} {\bf C} ,  
\] 
where ${\bf S} =( S_{ef} )_{e,f \in D(G)} $ is given by 
\[
S_{ef} =\left\{
\begin{array}{ll}
1 & \mbox{if $f= e^{-1} $, } \\
0 & \mbox{otherwise. }
\end{array}
\right. 
\]
The matrix ${\bf S} $ is called the {\em shift operator}.

The time evolution of a quantum walk on $G$ through ${\bf U} $ is given by 
\[
\psi {}_{t+1} = {\bf U} \psi {}_t . 
\]
Here, $\psi {}_{t+1}, \psi {}_t $ are the states. 
Note that the state $\psi {}_t $ is written with respect to the initial state $\psi {}_0 $ as follows: 
\[ 
\psi {}_{t} = {\bf U}^t \psi {}_0 . 
\]
A quantum walk on $G$ with ${\bf U} $ as a time evolution matrix is called a 
{\em coined quantum walk} on $G$.

If $\alpha {}_u (e)= \frac{1}{\sqrt{d_u}} $ for $e \in D(u)$, then 
the time evolution matrix ${\bf U} $ is called the {\em Grover matrix} of $G$, and 
a quantum walk on $G$ with the Grover matrix as a time evolution matrix is called a 
{\em Grover walk} on $G$.  
Thus, the {\em Grover matrix} ${\bf U} ={\bf U} (G)=( U_{ef} )_{e,f \in D(G)} $ 
of $G$ is defined by 
\[
U_{ef} =\left\{
\begin{array}{ll}
2/d_{t(f)} (=2/d_{o(e)} ) & \mbox{if $t(f)=o(e)$ and $f \neq e^{-1} $, } \\
2/d_{t(f)} -1 & \mbox{if $f= e^{-1} $, } \\
0 & \mbox{otherwise}
\end{array}
\right. 
\]

Let $G$ be a connected graph with $n$ vertices and $m$ edges. 
Then the $n \times n$ matrix ${\bf T}_n (G)=( T_{uv} )_{u,v \in V(G)}$ is given as follows: 
\[
T_{uv} =\left\{
\begin{array}{ll}
1/( \deg {}_G u)  & \mbox{if $(u,v) \in D(G)$, } \\
0 & \mbox{otherwise.}
\end{array}
\right.
\] 
Note that the matrix ${\bf T} (G)$ is the transition probability matrix of the simple random walk on $G$.

\begin{theorem}[Konno and Sato]
Let $G$ be a connected graph with $n$ vertices $v_1 , \ldots , v_n $ and $m$ edges. 
Then the characteristic polynomial for the Grover matrix ${\bf U}$ of $G$ is given by 
\[
\begin{array}{rcl}
\det ( \lambda {\bf I}_{2m} - {\bf U} ) & = & ( \lambda {}^2 -1)^{m-n} \det (( \lambda {}^2 +1) {\bf I}_n -2 \lambda {\bf T} (G)) \\
\  &   &                \\ 
\  & = & \frac{( \lambda {}^2 -1)^{m-n} \det (( \lambda {}^2 +1) {\bf D} -2 \lambda {\bf A} (G))}{d_{v_1} \cdots d_{v_n }} .   
\end{array}
\]
\end{theorem}

From Theorem 3, the following equation for the Grover matrix on a graph is obtained.

\newtheorem{corollary}{Corollary} 
\begin{corollary}
Let $G$ be a connected graph with $n$ vertices $v_1 , \ldots , v_n $ and $m$ edges. 
Then the characteristic polynomial for the Grover matrix ${\bf U}$ of $G$ is given by 
\[
\begin{array}{rcl}
\det ( {\bf I}_{2m} -u {\bf U} ) & = & (1- u^2 )^{m-n} \det ((1+ u^2 ) {\bf I}_n -2u {\bf T} (G)) \\
\  &   &                \\ 
\  & = & \frac{(1- u^2 )^{m-n} \det ((1+ u^2 ) {\bf D} -2u {\bf A} (G))}{d_{v_1} \cdots d_{v_n }} .   
\end{array}
\]
\end{corollary}

\section{The generalized Ihara zeta function of $\mathbb{Z}$} 

Let the {\em cycle graph} $C_n $ be the connected 2-regular graph with $n$ vertcies. 
If $n \rightarrow \infty$, then the limit of $C_n $ is the one-dimensional integer lattice $\mathbb{Z} $. 
Then we consider the Grover walk on $\mathbb{Z} $. 
This quantum walk is a 
\[
\left[
\begin{array}{cc}
1 & 0 \\
0 & 1 
\end{array}
\right]
\]
free quantum walk.

Let ${\bf U}^{(s)}_n $  be the Grover matrix on $C_n $ and ${\bf P}^{(s)}_n $ the transition probability matrix of the simple random walk on $C_n $. 
Then we have 
\[
{\ P}^{(s)}_n = \frac{1}{2} {\bf A} (C_n ) . 
\]
By Corollary 1, we have 
\[
\det ({\bf I}_{2n} -u {\bf U}^{(s)}_n ) =(1- u^2 )^{n-n} \det ((1+u^2 ) {\bf I}_n -2u {\bf P}^{(s)}_n )
=(1- u^2 )^{n-n} \det ( {\bf I}_n -u {\bf A} (C_n )+ u^2 {\bf I}_n) . 
\]
That is, 
\[
{\bf Z} ( C_n , u)^{-1} = \det ({\bf I}_{2n} -u {\bf U}^{(s)}_n ) . 
\]

By the fact that, for a finite graph, the classical Ihara zeta function is just the above
Ihara zeta function raised to the power equaling the number of vertices, we have 
\[
\begin{array}{rcl}
\zeta {}_{C_n } (u)^{-1} & = & {\bf Z} (C_n , u)^{-1/n} = \det ( {\bf I}_{2m} -u {\bf U}^{(s)}_n  )^{1/n} \\ 
\  &   &                \\ 
\  & = & \{ \det ((1+u^2 ) {\bf I}_n -2u {\bf P}^{(s)}_n ) \} {}^{1/n } \\ 
\  &   &                \\ 
\  & = & \{ \prod_{ \lambda \in {\rm Spec}( {\bf P}^{(s)}_n )} ((1+u^2 )-2u \lambda ) \} {}^{1/n} \\ 
\  &   &                \\ 
\  & = & \exp [ \log \{ ( \prod_{ \lambda \in {\rm Spec}( {\bf P}^{(s)}_n )} ((1+u^2 )-2u \lambda ))^{1/n} \} ] \\ 
\  &   &                \\ 
\  & = & \exp [ \frac{1}{n} \sum_{\lambda \in {\rm Spec} ( {\bf P}^{(s)}_n )} \log ((1+u^2 )-2u \lambda ) ] . 
\end{array}
\]

Now, since 
\[ 
\det ({\bf I}_{2n} -u {\bf U}^{(s)}_n )= {\bf Z} (C_n , u)^{-1} = (1- u^n)^2 , 
\]
we have 
\[
\det ( \lambda {\bf I}_{2n} - {\bf U}^{(s)}_n )=( \lambda {}^n -1)^2 .  
\]
Thus,  
\[
{\rm Spec} ( {\bf U}^{(s)}_n )= \{ [ e^{i \theta {}_0 } ]^2 , [ e^{i \theta {}_1 } ]^2 , \ldots , [ e^{i \theta {}_{n-1} } ]^2 \} , 
\]
where ${\rm Spec} ({\bf F})$ is the spectra of a square matrix ${\bf F}$, and 
\[
\theta {}_k = \frac{2 \pi k}{n} \ (k=0,1, \ldots , n-1) . 
\]
Furthermore, by Theorem 3 (Konno-Sato Theorem), we obtain the following spectral mapping theorem: 
\[
\det ( \lambda {\bf I}_{2n} - {\bf U}^{(s)}_n )=(2 \lambda )^n ( \lambda {}^2 -1)^0 \det ( \frac{ \lambda +1}{2 \lambda } {\bf I}_n - {\bf P}^{(s)}_n ) 
=(2 \lambda )^n \det ( \frac{ \lambda + \overline{ \lambda } }{2} {\bf I}_n - {\bf P}^{(s)}_n ) . 
\]
If $ \lambda =e^{i \theta {}_k } ( k=0,1, \ldots , n-1)$, then we obtain  
\[
\frac{ \lambda + \overline{ \lambda } }{2} = \cos \theta {}_k . 
\]
That is,  
\[
{\rm Spec} ( {\bf P}^{(s)}_n )= \{ [ \cos \theta {}_k ]^1 \mid k=0,1, \ldots , n-1\} .  
\]
Thus, we have 
\[
\zeta {}_{C_n } (u)^{-1} = \exp \Bigl[ \sum^{n-1}_{k=0} \log \Bigl( (1+u^2 )
-2u \cos \Bigl( \frac{2 \pi k}{n} \Bigr) \Bigr) \frac{1}{2 \pi }\times \frac{2 \pi }{n} \Bigr] . 
\]
When $n \rightarrow \infty$, then we have 
\[ \displaystyle
\begin{array}{rcl}
\displaystyle \lim {}_{n \rightarrow \infty} \zeta {}_{C_n } (u)^{-1} & = & \displaystyle 
\exp \Bigl[ \int^{2 \pi}_0 \log \Bigl( (1+u^2 )
-2u \cos x \Bigr) \frac{dx}{2 \pi } \Bigr] \\ 
\  &   &                \\ 
\  & = & \displaystyle \frac{u^2 +1+| u^2 -1|}{2} = \left\{
\begin{array}{ll}
1 & \mbox{if $|u| <1$, } \\
u^2  & \mbox{if $|u| \geq 1$. }
\end{array}
\right. 
\end{array}
\]

On the other hand, we see that 
\[
\zeta {}_{\mathbb{Z}} (u)=1 ,  
\] 
since $\mathbb{Z} $ has no reduced cycle. 
Therefore we obtain the following result.

\begin{theorem}  
\[
\lim {}_{n \rightarrow \infty } \zeta {}_{C_n } (u)=\zeta {}_{\mathbb{Z}} (u) \ for \  |u| <1.   
\]
\end{theorem}

From now on we consider a relation between Theorem 1.3 given by Chinta et al. [2] and Theorem 3. 

Their result for $ \mathbb{Z} $ case gives 
\[ \displaystyle 
\begin{array}{rcl}
\displaystyle  \zeta {}_{\mathbb{Z}} (u)^{-1} & = & \displaystyle  (1-u^2 )^{(1-1)/2} \exp ( \int \log (1-2u+ u^2 +u \lambda ) d \mu ( \lambda )) \\ 
\  &   &                \\ 
\  & = & \displaystyle \exp ( \int \log (1-2u  +u^2 + \lambda u ) d \mu ( \lambda )) \  \  \  \  \ \  (1) . 
\end{array} 
\]
Noting that $\Delta (C_n )= {\bf D} - {\bf A} (C_n )=2( {\bf I}_n - {\bf P}^{(s)}_n ) $, 
we have 
\[
\begin{array}{rcl}
\zeta {}_{C_n } (u)^{-1} & = & {\bf Z} (C_n , u)^{-1/n} = \det ( {\bf I}_{2m} -u {\bf U}^{(s)}_n  )^{1/n} \\ 
\  &   &                \\ 
\  & = & \{ \det ((1-2u+u^2 ) {\bf I}_n +u \Delta (C_n )) \} {}^{1/n } . 
\end{array}
\] 
Thus a similar argument in the proof of Theorem implies 
\[
\lim {}_{n \rightarrow \infty } \zeta {}_{C_n } (u)^{-1} = 
\exp \Bigl[ \int^{2 \pi }_0 \log ((1-2u+u^2 )+2u (1-\cos x )) \frac{dx}{2 \pi } \Bigr] . 
\]
Remark that the right-hand side of this equality is noting but that of Eq. (1). 
Then we have

\begin{corollary}
\[  
\displaystyle 
\lim {}_{n \rightarrow \infty } \zeta {}_{C_n } (u)^{-1} = \zeta {}_{\mathbb{Z}} (u)^{-1} 
= \exp \Bigl( \int \log (1-2u  +u^2 + \lambda u ) d \mu ( \lambda ) \Bigr) , 
\] 
where 
\[
\lambda d \mu ( \lambda ) \sim 2(1 -\cos x) \frac{dx}{2 \pi }  \  \  \  \  \  on \ [0, 2 \pi ) . 
\]
\end{corollary}

\end{document}